\documentclass[12pt]{amsart}
\usepackage{diagbox}
\usepackage{mathtools}
\usepackage{bbm}
\usepackage{amsmath}
\allowdisplaybreaks

\mathtoolsset{showonlyrefs}

\renewcommand{\emph}[1]{\textit{#1}}
\usepackage{enumerate,amsmath,amsthm,latexsym,amssymb}
\usepackage{color}\usepackage{graphicx}

\definecolor{brown}{cmyk}{0, 0.72, 1, 0.45}
\definecolor{grey}{gray}{0.5}

\def\bu{{\bf u}}
\def\bv{{\bf v}}

\newcommand{\old}[1]{}

\newcounter{rot}%\addtocounter{rot}{1}, \therot

\newcommand{\card}[1]{\left|#1\right|}

\newcommand{\ignore}[1]{}

\def\bM{{\bf M}}

\def\cA{{\mathcal A}}
\def\cB{{\mathcal B}}
\def\cS{{\mathcal S}}

\newcommand{\set}[1]{\left\{#1\right\}}

\newcommand{\proofstart}{{\noindent \bf Proof\hspace{2em}}}
\newcommand{\proofend}{\hspace*{\fill}\mbox{$\Box$}\\ \medskip\\ \medskip}

\def\ii_(#1,#2){i_{#1}^{#2}}

\def\bx{{\bf x}}

\def\cM{\mathcal{M}}
\def\a{\alpha}
\def\b{\beta}
\def\d{\delta}

\def\e{\varepsilon}
\def\f{\phi}

\def\g{\gamma}

\def\z{\zeta}

\def\th{\theta}

\def\l{\lambda}
\def\m{q}
\def\n{p}

\def\r{\rho}

\def\s{\sigma}

\def\om{\omega}

\def\be{{\bf e}}

\def\bos{{\boldsymbol \eta}}

\newcommand{\rdup}[1]{\left\lceil #1 \right\rceil}

\def\bss{{\boldsymbol \sigma}}
\def\bsa{{\boldsymbol \alpha}}
\def\bsb{{\boldsymbol \beta}}
\def\brho{{\boldsymbol \rho}}
\def\cE{\mathcal{E}}

\def\cI{{\mathcal{I}}}

\newcommand{\brac}[1]{\left( #1 \right)}

\def\E{{\bf E}}

\renewcommand{\Pr}{\operatorname{\bf Pr}}
\newcommand\bfrac[2]{\left(\frac{#1}{#2}\right)}

\def\bC{{\bf C}}

\def\bm{{\bf m}}
\newtheorem{theorem}{Theorem}[section]

\newtheorem{lemma}[theorem]{Lemma}

\newtheorem{remthm}[theorem]{Remark}

\newenvironment{remark}{\begin{remthm}}{\end{remthm}}%
\newcounter{thmtemp}

\newcommand{\nospace}[1]{}

\def\path{\operatorname{PATH}}

\newcommand{\beq}[2]{\begin{equation}\label{#1}#2\end{equation}}

\def\leb{\leq_b}
\def\geb{\geq_b}
\def\bA{{\bf A}}
\def\bD{\bf D}

\newcommand{\mult}[2]{\begin{multline}\label{#1}#2\end{multline}}
\def\bX{{\bf X}}
\newcommand{\sbs}{\subset}
\newcommand{\stm}{\setminus}

\def\bB{{\bf B}}
\def\bc{{\bf c}}

\def\bI{{\bf I}}
\def\br{{\bf r}}
\def\bL{{\bf L}}
\def\ba{{\bf a}}
\def\bR{{\bf R}}
\def\bC{{\bf C}}

\def\bb{{\bf b}}
\def\bd{{\bf d}}
\def\GF{\mathbf{GF}}
\def\bt{{\boldsymbol \theta}}
\def\bhA{\widehat{\bA}}
\begin{document}
\title{Minors of a random binary matroid}
\author{Colin Cooper}\thanks{Research supported in part by EPSRC grant EP/M005038/1}
\author{Alan Frieze}\thanks{Research supported in part by NSF Grants DMS1362785, CCF1522984 and a grant(333329) from the Simons Foundation}
\author{Wesley Pegden}\thanks{Research supported in part by NSF grant
    DMS1363136}
\date{\today}
\begin{abstract}
Let ${\bf A}$ be an $n\times m$ matrix over $\mathbf{GF}_2$ where each column consists of $k$ ones, and let $M$ be an arbitrary fixed binary matroid. The matroid growth rate theorem implies that there is a constant $C_M$ such that $m\geq C_M n^2$ implies that the binary matroid induced by {\bf A} contains $M$ as a minor. We prove that if the columns of ${\bf A}={\bf A}_{n,m,k}$  are chosen \emph{randomly}, then there are constants $k_M, L_M$ such that $k\geq k_M$ and $m\geq L_M n$ implies that ${\bf A}$ contains $M$ as a minor w.h.p.
\end{abstract}
\maketitle
\section{Introduction}
There is by now a vast and growing literature on the asymptotic properties of random combinatorial structures. First and foremost in this context are Random Graphs and Hypergraphs, see \cite{B}, \cite{FK} and \cite{JLR} for books on this subject. Random groups in their own right and in the guise of random permutations are included in this. Going further afield into Algebraic Geometry we see a recent surge of interest in Random Simplicial Complexes, initiated by the paper of Linial and Meshulam \cite{LM06}. See Kahle \cite{K15} for a recent survey. Another area of interest in this vein is that of Random Matroids. This paper concerns one aspect of these. For the basic facts on matroids see Welsh \cite{Welsh} or Oxley \cite{Ox}.  Basically we see that two models of a random matroid have been considered so far. 

In the first model a matroid is chosen uniformly at random from the set of all matroids with $n$ elements, see for example Oxley, Semple, Warshauer and Welsh \cite{OSWW}. Recently, there have been some breakthrough results in this subject. Bansal, Pendavingh and van der Pol \cite{BPP} give a very close estimate for $\log\log m_n$ where $m_n$ is the number of matroids on a fixed ground set with $n$ elements. And Nelson \cite{Nel} showed that almost all matroids are non-representable. Pendavingh and van der Pol \cite{PP} considered random matroids of rank $r$ and showed that almost all $r$-sets will be bases in this model.

The second model considers representable matroids. Given a matrix \bA\ we let $\cM(\bA)$ denote the representable matroid with ground set equal to the columns of \bA\ and independence given by linear independence.  An example of this model is the space of $n\times m$ matrix with entries chosen independently and uniformly from $\GF_q$, see for example Kelley and Oxley \cite{KO}. 

The random graph $G_{n,m}$ can be identified with a random $n\times m$ (0,1)-matrix $\bA_{n,m,2}$ where each row represents a vertex and each column has exactly two ones and defines an edge. If the entries are considered to be in $\GF_2$ and the ones in each column are chosen at random, then we have a matrix representation of a random graph and a random graphic matroid. 

The columns of $\bA_{n,m,2}$ define a (random) graphic matroid.  If we want to generalize this to random sample from a larger class of binary matroids, then one natural way is to take  $k$ random ones instead of 2 ones in each column, to obtain the random matrix $\bA_{n,m,k}$, the vertex-edge incidence matrix of a random $k$-uniform hypergraph. It is this model of a random binary matroid that is the subject of this paper. 

Many properties of a matroid are determined by whether or not it contains some particular fixed matroid as a minor. For example a binary matroid is regular if and only if it does not contain the Fano plane or its dual as a minor, see Tutte \cite{Tu}. We are interested in the event that $\bA_{n,m,k}$ contains a fixed binary matroid $M$ as a minor.  The matroid growth rate theorem of Geelen, Kung and Whittle \cite{growth} implies that there is a constant $C_M$ depending only on $M$ such that \emph{any} binary matroid of rank $n$ on $m>C_M n^2$ elements must contain $M$ as a minor (see also Kung  \cite{Kung}.  We prove that (when $k$ is large), for the random matroid induced by $\bA_{n,m,k}$ this quadratic condition can be replaced by a linear one. We prove the following:
\begin{theorem}\label{th1}
Let $M$ be a fixed binary matroid. Then there exist constants $k_M,L_M$ such that if $k\geq k_M$ and $m\geq L_Mn$ then w.h.p.\footnote{A sequence of properties $\cE_n,n\geq 1$ is said to hold {\em with high probability} (w.h.p.) if $\lim_{n\to\infty}\Pr(\cE_n)=1$.} $\cM_{n,m,k}$ contains $M$ as a minor.
\end{theorem}

We briefly recall the definition of the minor relation for matroids.  Given a matroid $\cM$ on the ground set $E$ and with the family $\cI$ of independent sets, for $X\subseteq E$, the deletion $\cM\stm X$ is the matroid on $E\stm X$ whose independent sets consist of $\{I\in \cI:\:\: I\sbs E\stm X\}$.  The contraction $\cM/X$ $(X\in \cI)$ is the matroid on $E\stm X$ whose independent sets are $\set{I\sbs E\stm X: \:\: I\cup X\in \cI}$.  $M$ is a minor of $\cM$ if it can be obtained from $\cM$ by deletion and contraction operations. (For $X\notin \cI$, the contraction can be defined by $\cM/X:=(\cM/Y)\stm X$, where $Y$ is any basis of $X$.)

Theorem \ref{th1} is related to the result of Altschuler and Yang \cite{AY}. They prove that if matrix \bM\ is an $n(m)\times m$ matrix with random entries in $\GF_q$ and $m-n(m)\to\infty$ then w.h.p. the matroid associated with \bM\ contains any fixed minor. This can be related to our theorem on taking $q=2$ and $k=n(m)/2$. (We have reversed the roles of $m,n$ from their statement.) However, the results of \cite{AY} rely heavily on the fact that pre-multiplying a uniform random matrix in this model by a non-singular matrix yields another uniform random matrix. Our model lacks this property. Furthermore, multiplying $\bA_{n,m,k}$ by a non-singular matrix will not fix this property. This is because whatever matrix we use as a pre-multiplier, we will only have a sample space of size at most $\binom{n}{k}$ for the resulting column set, as opposed to $2^n$.

\section{Proof of Theorem \ref{th1}}
\subsection{Outline of our proof}\label{outline}
Fix $k$ and let the matrix $\bA_m=\bA_{n,m,k}$ have columns\\ $[\ba_1,\ba_2,\ldots,\ba_m]$ where $m=Kn$ for $K$ sufficiently large. Let $M$ be a fixed binary matroid and let $\bR_M=[\bm_1,\bm_2,\ldots,\bm_\m]$ be a representation of $M$ by a $\n\times\m$ matrix. Assume without loss of generality that $\bR_M$ has full row-rank $\n$. In this outline we will assume that $k$ is odd. There are some minor adjustments needed for $k$ even.

Let 
\beq{m1n1}{
n_1=n\text{ and }m_1=\frac{n}4. 
}
Denote the $n_1\times m_1$ matrix consisting of the first $m_1$ columns of $\bA_m$ by $\bX$. It follows from Theorem 1 of Cooper \cite{Co1} that w.h.p. the columns of $\bX$ are linearly independent.

We will use results on hypergraph cores to find a sub-matrix $\bB_1$ of \bX\ that has $n_2$ rows and $m_2$ linearly independent columns where $n_2$ is close to $n_1$ and $m_2\ge n/5$, and with the property that $\bB_1$ has $k$ random ones in each column and at least $k/10$ ones in each row. 

We extend $\bB_1$ to an $n_2\times n_2$ non-singular submatrix $\bB$ of $\bA_m$ which again has exactly $k$ ones in each column, as follows. Let $I_1$ denote the index set of the rows of $\bB_1$. We extend $\bB_1$ by choosing $Ln$ columns of $\bA_m$, disjoint from \bX\ to create a submatrix \bL. Here $L$ is a sufficiently large constant. These columns will only have ones in rows indexed by $I_1$. Again using properties of hypergraph cores, we show that w.h.p. \bL\ contains a submatrix $\bL_1$ which has $n_3$ rows and $m_3$ columns which (i) has full row rank and (ii) each row has at least $\z kL$ ones. Here $n_3$ is close to $n_2$ and $0<\z<1$. We next obtain $\bL_3$ from $\bL_1$ by adding $n_2-n_3$ rows of zeros. We then argue that the matrix $\bL_2=[\bB_1:\bL_3]$ has $n_2$ rows and has full row rank w.h.p. The matrix \bB\ is an arbitrary extension of $\bB_1$ to a square non-singular $n_2\times n_2$ sub-matrix of $\bL_2$.

We then argue that w.h.p. the rows of $\bB^{-1}$ have between $\e_0n_2=\frac{1}{2}e^{-k}n_2$ and $n_2-\e_0n_2$ ones. 

We let $\bhA$ be the $n_2\times m_3$ submatrix of $\bA_m$ whose rows are the rows of $\bB$, and whose columns are those columns of $\bA_m$ which have ones only in rows of $\bB$. Note that $\cM(\bhA)$ is a minor of $\cM(\bA_m)$.  Now write $\bhA=[\bB:\bM]$ and consider the matrix $\bhA_1=[\bI:\bM_1]$ for $\bM_1=\bB^{-1}\bM$, where we assume that the first $n_2$ columns form the $n_2\times n_2$ identity matrix. Suppose that $\bM_1$ contains a submatrix equal to our target matrix $\bR_M$. Then we are done. Indeed, suppose w.l.o.g. that $\bR_M$ lies in the first $\n$ rows and the first $\m$ columns of $\bM_1$. Then we get $M$ as a minor of $\cM(\bhA)$ (and hence of $\cM(\bA_m)$) by deleting the first $\n$ columns of $\bB$ and the last $m_3-n_2-\m$ columns of $\bM$ and contracting the last $n_2-\n$ columns of $\bB$, as we explain next.

Recall that a minor of $\bA_m$ is obtained by deleting and contracting columns.  Recall from the definition of contraction that if $S$ denotes an independent set (of column indices), then a set $T$ (of column indices) disjoint from $S$ is independent in the contraction $\cM/S$ iff $S\cup T$ is an independent set (of columns) in $\cM$. 

Contraction is simple if the columns $S$ are a subset of the columns of an identity matrix $\bI=\bI_{n_2}$. In view of this, we pre-multiply $\bhA=[\bB:\bM]$ by $\bB^{-1}$ to obtain $\bhA_1=[\bI:\bM_1]$. Pre-multiplying by a non-singular matrix does not change the underlying matroid, seeing as column dependence/independence is preserved. We can assume that the first $n_2$ columns form the $n_2\times n_2$ identity matrix $\bI$. If we contract a set $S$ of the columns of $\bI$, then a representation of the contracted matroid is given by deleting the $|S|$ rows of $\bhA_1$ that have a one in a column of $S$ to obtain a matrix $\bhA_2$. In which case we see that a set $T$ of columns of $\bhA_2$ is independent in $\bhA_2$ if and only if the set of columns corresponding to $S\cup T$ is independent in the matroid represented by $\bA_m$.

To prove that $\bR_M=[\bm_1,\bm_2,\ldots,\bm_\m]$ appears as a submatrix of $\bM_1$, we will consider $\bB^{-1}\bc$ where $\bc$ is a random column of $\bhA$ outside of the $n/4+Ln$ columns considered so far in the construction of \bB. For a set $R$ of rows and a column \bx\ of $\bB^{-1}\bhA$, let $\f_R(\bx)$ be the column \bx\ restricted to the rows $R$. We argue next that we can find $R$ of size $p$ such that 
\beq{done1}{
\Pr(\f_R(\bB^{-1}\bc)=\bm_j)=\Omega(1)\text{ for }1\leq j\leq \m.
}
This means that w.h.p. we can find a copy of each column of $\bR_M$ by searching through $\om$ random columns, where $\om=o(n)$ is any function tending to infinity with $n$. 

To justify \eqref{done1}, let $S_i$ denote the support of the $i$th row of $\bB^{-1}$. Our strategy for analyzing $\bB^{-1}\bc$ is to show that there is a set $R$ of $\n$ rows of $\bB^{-1}$ and a partition $A_0,A_1,\ldots,A_\ell$ of $[n_2]$ such that for all $i\in R,1\leq j\leq \ell$, $S_i$ contains $A_j$ or is disjoint from it. There will be a corresponding $\n\times \ell,\ (0,1)$-matrix ${\bD}=({\bD}[i,j])$ with the following properties.  ${\bD}$ has full row rank and for some constants $0<\e_1\ll \e_0\ll 1$, ${\bD}[i,j]=1$ implies (i) $r_{i,k}=1$ for $k\in A_j,i\in R$, ( $\br_i=(r_{i,\cdot})$ being the $i$th row of $\bB^{-1}$), (ii) $|A_j|\geq \e_1n_2$ for $j\geq 0$.

Given $R$ and ${\bD}$ we proceed as follows: Let $\bc=(c_1,c_2,\ldots,c_{n_2})$ be a random column with $k$ 1's. Let $\bv$ satisfy ${\bD}\bv=\bm_1$ (the first column of $\bR_M$) and $v_j=0$ if $|A_j|<\e_1n$. We can assume that $\bv$ has at most $\n$ ones and that $k\geq\n$. Equation \eqref{done1} follows from  
\beq{111}{
\Pr(\f_R(\bB^{-1}\bc)=\bm_1)\geq\Pr(\bc_R=\bv)=\Omega(1),
}
where $\bc_R=(d_0,d_1,\ldots,d_\ell)$ and where $d_j=\sum_{l\in A_j}c_l$. 

The condition in \eqref{111} will be satisfied if exactly one element is chosen from each $A_j$ such that $v_j=1$ and the rest are chosen from $A_0=[n_2]\setminus \bigcup_{i\in R}S_i$. This has probability $\Omega(1)$. 

We will show in Section \ref{notenough} how to choose the set of $p$ rows $R$ so that they contain at least $\e_1n_2$ common zeros. Then in Section \ref{largeM} we will show that if $S_i^1=S_1,S_i^0=\bar{S}_i$ then the partition $A_{{\boldsymbol \sigma}}= \bigcap_{j=1}^pS_j^{\xi_j},{\boldsymbol \sigma}=(\xi_1,\xi_2,\ldots,\xi_\n)$ (as $(\xi_1,\xi_2,\ldots,\xi_\n)$ runs over $\set{0,1}^\n$) suffices as a partition. We will take $D[i,{\boldsymbol \sigma}]=1$ only if $\xi_i=1$ and $|A_{{\boldsymbol \sigma}}|\geq \e_1n_2$. 
\subsection{Some Notation}
We summarize here the meaning of some parameters. The reader might find this useful to refer back to.
\begin{enumerate}[(i)]
\item $\bB_1$ is the $n_2\times m_2$ submatrix derived from the first $m_1$ columns of $\bA_{n,m,k}$. Every column has $k$ ones and every row has at least $k/10$ ones. The columns of $\bB_1$ are linearly independent and the values $n_2,m_2$ satisfy \eqref{C1==}, \eqref{m2==} below. The set $I_1$ is the index set of rows of $\bB_1$.
\item \bL\ is an $n_2\times Ln$ submatrix of $\bA_{n,m,k}$, whose columns are disjoint from those of $\bB_1$.
\item $\bL_1$ is an $n_3\times m_3$ submatrix of $\bL$ which has rank $n_3$, where $n_3,m_3$ satisfy \eqref{C1==0}, \eqref{m2===}.  The rows of $\bL_1$ have index $I_2\subseteq I_1$.
\item $\bL_2$ is an $n_2$ row matrix that contains $\bB_1$ as a sub-matrix and has rank $n_2$ and many more than $n_2$ columns. It is therefore possible to find an $n_2\times n_2$ non-singular submatrix \bB\, that contains $\bB_1$ and is contained in $\bL_2$.
\item In general bold named variables are either matrices or vectors.
\end{enumerate}
We now give a detailed proof of Theorem \ref{th1}.
\subsection{Building $\bB_1$}
Consider the $k$-uniform hypergraph $H_1$ induced by the first $m_1=n/4$ columns $\bX$ of $\bA_{n,m,k}$. I.e. the hypergraph with a vertex for each row and where each edge $e_j,j\leq n/4$ corresponds to the column $\bc_j$ of \bX\ via $e_j$ contains an element $i\in [n]$ if and only if $\bX[i,j]=1$. $H_1$ is distributed as a random $k$-uniform hypergraph with $n_1$ vertices and $m_1$ edges. We show next that w.h.p. the $k/10$-core $C_1$ of $H_1$ is large. The $r$-core of a hypergraph $H=(V,E)$ is the largest set $S\subseteq V$ such that each $s\in S$ has degree at least $r$ in the sub-hypergraph of $H$ induced by $S$ i.e. each $s\in S$ lies in at least $r$ edges $e$, $e\subseteq S$. The $k/10$-core will provide us with a matrix $\bB_1$ with at least $k/10$ ones in each row.

We use some results on the cores of  random $k$-uniform hypergraphs (see e.g. Cooper \cite{CDCCores} or Molloy \cite{Mo}). Let $c=km_1/n_1=k/4$, and
let $x$ be the greatest solution to
\beq{x=}{
c=\frac{k}{4}=\frac{x}{\brac{1-e^{-x}\sum_{i=0}^{k/10-2}\frac{x^i}{i!}}^{k-1}}.
}
We will use a simple continuity argument to prove the existence of $x$ and bound it as in \eqref{x==} below.

It is known that w.h.p.,
\beq{C1=}{
n_2=|V(C_1)|\approx n_1\brac{1-e^{-x}\sum_{i=0}^{k/10-1}\frac{x^i}{i!}},
}
and
\beq{edgescore}{
m_2=|E(C_1)|\approx m_1\bfrac{x}{c}^{k/(k-1)}.
}
Here, $A(x)\approx B(x)$ stands for $A(x)=(1+o(1))B(x)$ as $x\to\infty$, $A(x)\gtrsim B(x)$ stands for $A(x)\geq(1+o(1))B(x)$ as $x\to\infty$.

We will first argue that for $k$ large we have
\beq{x==}{
\frac{k}{5}< x\leq \frac{k}{4}.
}
The upper bound follows directly from the definition \eqref{x=}. To prove the lower bound let
$$S(x)=\frac{k}{4}-
\frac{x}{\brac{1-e^{-x}\sum_{i=0}^{k/10-2}\frac{x^i}{i!}}^{k-1}}.$$
If $x \ge 2(i+1)$, then  $\frac{x^i}{i!} \le \frac{x^{i+1}}{2(i+1)!}$. 
Thus for $x \ge k/5$ and $\th=1,2$,
\[
\sum_{i=0}^{k/10-\th}\frac{x^i}{i!} \le  \frac{x^{k/10}}{(k/10)!}\le \bfrac{10xe}{k}^{k/10}
\]
and so
$$e^{-x}\sum_{i=0}^{k/10-\th}\frac{x^i}{i!} \le \brac{\frac{10x e}{k}\;  e^{-10 x/k}}^{k/10} \le 
\bfrac{2}{e}^{k/10}\quad\text{since }x\geq k/5.$$
Thus $S(k/5)>0$ for $k$ large. As  $S(k/4)<0$, the lower bound in \eqref{x==}  follows from the continuity of $S(x)$.
It then follows from \eqref{C1=} that w.h.p.
\beq{C1==}{
n_1\geq n_2=|V(C_1)|\ge n_1\brac{1-\frac{1}{k}}.
}
Similarly, using \eqref{edgescore} along with $c=k/4$ and $x> k/5$ from \eqref{x==}
gives
\beq{m2==}{
\frac{n}{4}\geq m_2 \gtrsim\frac{n_1}{4}
\bfrac{4}{5}^{k/(k-1)} \ge \frac{n_2}{5},
}
for $k$ large.

Now consider the submatrix $\bB_1$ of $\bX$ comprised of the columns corresponding to the edges of $H_1$ that are contained in $C_1$. The distribution of ones in $\bB_1$ is that each of the $m_2$ columns chooses $k$ random ones from $n_2$ rows, subject only to each row having at least $k/10$ ones. This is an interpretation of a standard result on cores of graphs being random subject to a lower bound on minimum degree. Let $I_1$ denote the index set of the rows of $\bB_1$. Thus $|I_1|=n_2$.
\subsection{Extending $\bB_1$ to a basis}\label{extendB}
We fix some sufficiently large constant $L>1$ and begin by choosing $Ln$ columns of $\bA_m$ disjoint from \bX\ to make a sub-matrix $\bL$. We choose the first $Ln$ columns following $\bX$ that have ones only in rows indexed by $I_1$. The probability that a random column only has ones in rows $I_1$ is $\frac{\binom{n_2}{k}}{\binom{n}{k}}=\Omega(1)$ and so w.h.p. we only need to examine $O(n)$ columns of $\bA_m$ in order to find these $Ln$ columns. Now let $0<\z<1$ be a small constant. Let now $H_2$ denote the $k$-uniform hypergraph induced by the columns of $\bL$ and let $C_2=C_2(H_2)$ denote its $\z Lk$-core.
Using \cite{CDCCores}, \cite{Mo} once again we see that we have to let $x$ be the greatest solution to
\beq{x=0}{
c=Lk=\frac{x}{\brac{1-e^{-x}\sum_{i=0}^{\z Lk-2}\frac{x^i}{i!}}^{k-1}}.
}
Then w.h.p.,
\beq{C1=0}{
n_3=|V(C_2)|\approx n_2\brac{1-e^{-x}\sum_{i=0}^{\z Lk-1}\frac{x^i}{i!}}.
}
We will next argue that for $k,L$ large we have
\beq{x==0}{
\frac{(1+\z)Lk}{2}\leq x\leq Lk.
}
The upper bound follows directly from the definition \eqref{x=0}. To prove the lower bound let now
$$S(x)=Lk-\frac{x}{\brac{1-e^{-x} \sum_{i=0}^{\z Lk-2}\frac{x^i}{i!}}^{k-1}}.$$
If $x\geq \frac{(1+\z)(i+1)}{2\z}$ then $\frac{x^i}{i!} \le \frac{\xi  x^{i+1}}{(i+1)!}$ where $\xi=\frac{2\z}{1+\z}<1$. Thus for $x\geq \frac{(1+\z)Lk}{2}$ and $\th=1,2$,
$$\sum_{i=0}^{\z Lk-\th}\frac{x^i}{i!}\leq \frac{1}{1-\xi}\cdot\frac{x^{\z Lk}}{(\z Lk)!}\leq \frac{1}{1-\xi}\bfrac{ex}{\z Lk}^{\z Lk},$$
and since $\eta=\frac{1+\z}{2\z}<1$, then
$$e^{-x}\sum_{i=0}^{\z Lk-\th}\frac{x^i}{i!} \le \frac{1}{1-\xi}\brac{\frac{ex}{\z Lk}e^{-x/(\z Lk)}}^{\z Lk}\le \frac{\brac{\eta e^{1-\eta}}^{\z Lk} }{1-\xi}. $$
Thus $S(\frac{(1+\z)Lk}{2})>0$ for large $k$ and the lower bound in \eqref{x==0} follows by continuity.

It then follows from \eqref{C1=0} that for large enough $L$, we have that w.h.p.
\beq{C1==0}{
n_2\geq n_3=|V(C_2)|\ge n_2\brac{1-e^{-2k}}.
}
Similarly, using a similar expression to \eqref{edgescore} along with $c=Lk$ and $x \ge (1+\z)Lk/2$ in \eqref{x==}
gives us that the number $m_3$ of edges in $C_2$ satisfies
\beq{m2===}{
Ln\geq m_3 \gtrsim L n_2
\bfrac{1+\z}{2}^{k/(k-1)} \text{ and so }m_3\ge \frac{4(1+\z)Ln_2}{9},
}
for $k$ large.

We argue next that w.h.p. the matrix $\bL_1$ induced by $C_2$ has rank $n_3$. For this we rely on the following lemma, which we will need for several purposes:
\begin{lemma}\label{lemrank}
Let $\bA=(\bA[i,j])$ be an $N\times M$ matrix over $\GF_2$ chosen  uniformly at random from matrices where each column has $k$ ones, and condition on the event that each row has at least $\g k\s$ ones, where $\g<1$ and $\g k>1$ and $\s=M/N=O(1)$. Let $\bsa$ be a fixed member of $\GF_2^M$. If $\cE_{s,\bsa}$ is the event that there exists a set $S$ of rows with $|S|=s$ whose sum is $\bsa$, then 
\begin{enumerate}[(a)]
\item
%If $1\leq s\leq Ne^{-k}$ then
\beq{zerosum}{
\Pr(\exists 1\leq s\leq Ne^{-k}:\cE_{s,\bsa})=O(N^{-K}).%\leb M^{1/2}\bfrac{se^{2k/3}}{N}^{\g k\s s/3}.
}
\item
If $\s\geq e^{5k}/(1-\g)^2$ then
\beq{zerosum1}{
\Pr(\exists Ne^{-k}<s<N:\cE_{s,\bsa})=O(N^{-K}).%O(e^{-\xi N})\qquad\text{ for some fixed }\xi>0.
}
\end{enumerate}
We can take $K=\g k\s/6$ in the above.
\end{lemma}
Thus $K$ can be made arbitrarily large, by taking $k$ sufficiently large. Note also that we exclude $|S|=N$ from the statement of the lemma, since this would be false with $\bsa$ equal to all ones ($k$ odd) or all zeros ($k$ even).

We apply the lemma to $\bL_1$ by taking $N=n_3,M=m_3$ where $\z=1/2$ and then let $\g$ be equal to $\frac{Ln_3}{2m_3}\in \left[\frac{k-1}{2k}(1-e^{-2k}),\frac{3}{4}\right]$. The bounds on $\g$ being justified by \eqref{m1n1}, \eqref{C1==}, \eqref{C1==0} and \eqref{m2===}. Assume that $L\ge10e^{5k}$, so that the lower bound on $\s$ in (b) is satisfied. We will now make the following:
\begin{description}
\item[Assumption A] $k$ is odd.
\end{description}
We will deal with the case of $k$ even in Section \ref{keven}. Now if $k$ is odd and $\bL_1$ does not have full row rank, then $\cE_{s,\bsa}$ occurs with $\bsa={\bf 0}$ for some $1\leq s\leq n_3$. But Lemma \ref{lemrank} implies that
\beq{rowrank}{
\Pr(\exists 1\leq s\leq n_3:\cE_{s,\bsa}\text{ occurs})=O(n^{-K}).
}
%For $k$ odd, we can take $s=n_3$ in \eqref{rowrank}, but not for $k$ even. For $k$ even we only claim that $\bL_1$ has rank $n_3-1$ w.h.p. We therefore only claim that Lemma \ref{lemrank} gives us \eqref{rowrank} for $s<n_3$.

So, w.h.p. we have found an $n_3\times m_3$ matrix $\bL_1$ of rank $n_3$. Now consider the matrix $\bL_2=\left[\bB_1:\bL_3\right]$. Here $\bL_3$ is obtained from $\bL_1$ by adding $n_2-n_3$ rows of zeros. We claim that w.h.p. $\bL_2$ has rank $n_2$. Let $I_2\subseteq I_1$ be the row indices of $\bL_1$. Let the rows of $\bL_2$ be $\ba_1,\ba_2,\ldots,\ba_{n_2}$ and suppose that there exists $J\subseteq I_1$ such that $\sum_{i\in J}\ba_i=0$. Then we have $J\cap I_2=\emptyset$ else $\bL_1$ does not have rank $n_3$. We have $J\subseteq I_1\setminus I_2$ and then \eqref{C1==0} implies that $|J|\leq ne^{-2k}$. But then we obtain a contradiction from Lemma \ref{lemrank}(a) applied to the rows of $\bB_1$. 

Because $\bL_2$ has full row rank, we can obtain \bB\ as an extension of $\bB_1$ to an $n_2\times n_2$ non-singular sub-matrix of $\bL_2$. After this we order the columns of \bB\ so that the columns of $\bB_1$ come first.
\subsection{Proof of Lemma \ref{lemrank}}\label{secrank}
We first deal with small $s$. Suppose that $1\leq s\leq Ne^{-k}$. If $T\subseteq S\subseteq [N], |S|=s$, let $\cE_{j,T,S}$ denote the event that column $j$ of \bA\ has ones in all of the rows $T$ and zero's in the rows $S\setminus T$. Then where $\bsa=(\a_1,\a_2,\ldots,\a_M)$,
\beq{cE0}{
\Pr(\cE_{s,\bsa})\leq \sum_{S\subseteq [N], |S|=s}\ \sum_{\substack{d_j=\a_j\text{mod } 2,j\in[M]\\ d_1+d_2+\cdots+d_{M}\geq \g k\s s}}\ \sum_{S_j\subseteq S, |S_j|=d_j} \Pr\brac{\bigcap_{j=1}^{M}\cE_{j,S_j,S}}.
}
{\bf Explanation:} We sum over sets $S$ and then for each $j\in [M]$ we fix the number of ones  $d_j=|\set{i\in S:\bA[i,j]=1}$ of column $j$ that appear in the rows $S$. We then choose the rows $S_j$ where these ones appear and multiply by the probability that things are just so.

To estimate the probabilities in the RHS of \eqref{cE0} we will use the following model: we choose \bX\ uniformly from $[N]^{kM}$. Then column $i$ of \bA\ contains a one in positions $X_{k(i-1),j}$ for $1\leq j\leq k$ and $1\leq i\leq M$. It is possible that $X_{i,j_1}=X_{i,j_2}$ for some $i,j_1,j_2$. Let $\cS$ be the event that this does not happen. Then
\beq{norepeats}{
\Pr(\cS)\geq \brac{1-\frac{\binom{k}{2}}{N}}^M\geq e^{-k^2\s}.
}
{\bf Explanation:} $\frac{\binom{k}{2}}{N}$ bounds from above the probability that a fixed column contains a repeat, by the expected number of repeats. Each column is independently generated and \eqref{norepeats} follows.

Thus $\Pr(\cS)=\Omega(1)$ and events involving \bX\ that occur w.h.p. will also occur w.h.p. if we condition on $\cS$.

We see next that given $\cS$, each matrix with exactly $k$ ones in each column is equally likely. Indeed, each such matrix arises from the same number $(k!)^M$ of choices of \bX. Thus we can use \bX\ to generate our matrix \bA\ in a uniform way. It remains to deal with the lower bounds on row sums.

The row-sums $\r_i=\card{\set{(a,b):X_{a,b}=i}|,1\leq i\leq N}$ will be independent Poisson random variables, subject to $\r_i\geq\g k\s,i\in[N]$ and $\r_1+\r_2+\ldots+\r_{N}=kM$. This was proved in \cite{AFP} where the lower bound of $\g k\s$ is replaced by 2. We include a proof in an appendix for completeness. Thus
\beq{rho=l}{
\Pr(\r=l)=\frac{\l^l}{l!f_{\g k\s}(\l)}\text{ where }f_a(\l)=e^\l-\sum_{i=0}^{a-1}\frac{\l^i}{i!}.
}
Here we choose $\l$ so that $\E(\r)=k\s$, which implies that
\beq{ratio}{
\frac{\l f_{\g k\s-1}(\l)}{f_{\g k\s}(\l)}=k\s.
}
This choice of $\l$ ensures that $\Pr(\r_1+\r_2+\ldots+\r_{N}=kM)=\Omega(M^{-1/2})$. This follows from a version of the local central limit theorem, proved in \cite{AFP}.

It follows that for large $k$, we have
\beq{lamsize}{
\frac{k\s}{2}\leq \l\leq k\s\text{ and }f_{\g k\s}(\l)\geq e^{\g k\s/2}.
}
The upper bound in \eqref{lamsize} follows from the fact that $f_{\g k\s-1}(\l)>f_{\g k\s}(\l)$. The lower bound follows from the fact that if $k$ is large, then the RHS of \eqref{ratio} is large and then $\l$ approaches $k\s$ which is large. This then implies that $f_{\g k\s-1}(\l)$ approaches $f_{\g k\s}(\l)$ as $k$ grows.

Suppose now that we condition on the row sums  $\r_1=\th_1,\r_2=\th_2,\ldots,\r_{N}=\th_{N}$. Fix $S$ and $S_j\subseteq S,|S_j|=d_j,j\in [M]$. Then if $d_1+d_2+\cdots+d_{M}=d$ then
\mult{uprob}{
\Pr\brac{\bigcap_{j=1}^{M}\cE_{j,S_j,S}}\leq \frac{(kM-d)!}{(kM)!}\prod_{j=1}^{M}\prod_{i\in S_j}(\th_ik)=\\
\frac{1}{M^d}\prod_{l=0}^{d-1}\brac{1-\frac{l}{kM}}^{-1}\prod_{i\in S}\th_i^{\th_i} \leq \frac{e^{d^2/kM}}{M^d}\prod_{i\in S}\th_i^{\th_i}.
}
{\bf Explanation of \eqref{uprob}:} 
The conditioned model involves a vector $\bX\in [N]^{kM}$ that can be viewed as a random permutation of $\r_i$ copies of $i$ for $i\in [N]$. We can assume that these copies are distinguishable. Then, if $i\in S_j$ and $(i_1,j_1),\ldots,(i_l,j_l)$ represent prior assignments,
\beq{explain17}{
\Pr(\bA[i,j]=1\mid \bA[i_1,j_1]=1,\ldots,\bA[i_l,j_l]=1)\leq \frac{k\th_i}{kM-l}.
}
To see \eqref{explain17}, observe that there are at most $k$ positions in \bX\ that give us $\bA[i,j]=1$ and for each there are at most $\r_i$ out of $kM-l$ equally likely choices of being $i$. The second term in equation \eqref{uprob} follows. 

Next let
$$D_\ell=\set{\bd=(d_1,d_2,\ldots,d_{M}):d_j=\a_j\,\text{mod } 2, d_j\leq k, j\in[M],\sum_{j\in [M]}d_j=\ell}$$
and
$$E_\ell=\set{\bt=(\th_i,i\in S):\sum_{i\in S}\th_i=\ell,\, \th_i\geq \g k\s,i\in S}.$$
Note that
\beq{ineqDE}{
|D_\ell|\leb \binom{M+\ell/2-1}{\ell/2-1}\leq \frac{M^{\ell/2}e^{\ell^2/4M}}{(\ell/2)!}
\text{ and }|E_\ell|=\binom{\ell-\g k\s s+s-1}{s-1}<2^\ell.
}
Here the notation $A\leb B$ is used in place of $A=O(B)$.

The first inequality in \eqref{ineqDE} is obtained as follows: Let $d_j'=(d_j-1)/2$ if $\a_j=1$ and let $d_j'=d_j/2$ if $\a_j=0$. Then $\sum_j d_j'=(\ell-\ell_1)/2$ where $\ell_1$ is the number of $\a_j$ equal to one. Knowing $\bsa$, which is fixed, we can re-construct the $d_j$'s from the $d_j'$'s. This explains the binomial coefficient. After this we use 
$$B!\binom{A+B}{B}=A^B\prod_{i=0}^{B-1}\brac{1+\frac{B-i}{A}}\leq A^Be^{B^2/A}.$$

Plugging \eqref{uprob} into \eqref{cE0} we obtain,
\begin{align}
&\Pr(\cE_{s,\bsa})\\
&\leq \sum_{S\subseteq [N], |S|=s}\sum_{\ell=\g k\s s}^{kM}\ \sum_{\bd\in D_\ell}\ \sum_{S_j\subseteq S, |S_j|=d_j}\ \sum_{\bt\in E_\ell}\Pr(\r_i=\th_i,i\in S)\times\frac{e^{\ell^2/kM}}{M^\ell}\prod_{i\in S}\th_i^{\th_i}\\
&\leb M^{1/2}\sum_{S\subseteq [N], |S|=s}\sum_{\ell=\g k\s s}^{kM}\ \sum_{\bd\in D_\ell}\ \sum_{S_j\subseteq S, |S_j|=d_j}\ \sum_{\bt\in E_\ell}\prod_{i\in S}\frac{\l^{\th_i}}{\th_i!f_{\g k\s}   (\l)}\frac{e^{\ell^2/kM}}{M^\ell}\prod_{i\in S}\th_i^{\th_i}.\\
&\leb M^{1/2}\binom{N}{s}\sum_{\ell=\g k\s s}^{kM}\ \sum_{\bd\in D_\ell}\ \sum_{\bt\in E_\ell}\frac{\l^\ell}{f_{\g k\s}(\l)^s}\bfrac{se^{\ell/kM}}{M}^\ell\prod_{i=1}^s\frac{\th_i^{\th_i}}{\th_i!}\\
\end{align}
The $M^{1/2}$ factor in the third line follows from our choice of $\l$ from \eqref{ratio}. To obtain the last line we used $\sum_{S_j\subseteq S, |S_j|=d_j}1\leq s^{d_1+\cdots+d_M}=s^\ell$.

Thus, 
\begin{align}
\Pr(\cE_{s,\bsa})&\leb M^{1/2}\binom{N}{s}\sum_{\ell=\g k\s s}^{kM}\ \sum_{\bd\in D_\ell}\ \sum_{\bt\in E_\ell}\frac{\l^\ell}{f_{\g k\s}(\l)^s}\bfrac{se^{1+\ell/kM}}{M}^\ell\\
&\leb M^{1/2}\bfrac{Ne}{s}^s\sum_{\ell=\g k\s s}^{kM}\frac{M^{\ell/2}e^{\ell^2/4M}}{(\ell/2)!} \frac{(2k\s)^\ell}{e^{\g k \s s/2}}\bfrac{se^{1+\ell/kM}}{M}^\ell\\
&\leb M^{1/2}\bfrac{Ne}{s}^s\sum_{\ell=\g k\s s}^{kM}\bfrac{e^{k/3}\s s}{\ell^{1/2}M^{1/2}}^\ell e^{-\g k\s s/2}, \label{SUM}
\end{align}
since $k$ is large. Now if $u_\ell$ is the summand in \eqref{SUM} then
$$\frac{u_{\ell}}{u_{\ell-2}}\leq \frac{e^{2k/3}\s^2s^2}{\ell M}\leq \frac{e^{2k/3}\s^2s^2}{\g ks\s^2N} =\frac{e^{2k/3}s}{\g kN}\leq \frac12,$$
since $\g k>1$.

Hence, since the largest term in the sum in \eqref{SUM} is at $\ell=\g k\s s$, it follows  that
\beq{abound}{
\Pr(\cE_{s,\bsa})\leb M^{3/2}\bfrac{Ne}{s}^s\bfrac{e^{2k/3}s}{N}^{\g k\s s/2} \leq  M^{3/2}\bfrac{se^{2k/3}}{N}^{\g k\s s/3}.
}
Summing the RHS of \eqref{abound} for $1\leq s\leq Ne^{-k}$ and taking $k$ large completes the proof of part (a) of the lemma.

Assume now that $Ne^{-k}\leq s\leq N/2$.  If the sum of the rows in $S$ is {\bf 0}, (resp. {\bf 1}), then no column has exactly one one (resp. exactly two ones) in the rows of $S$. Let these events be $\cA_{S,i},i=0,1$. If the ones in each column were generated completely at random then, with the aid of the Vandermonde identity,
\begin{align}
\Pr(\cA_{S,0})&=\brac{\sum_{i\neq 1}^{k}\frac{\binom{s}{i}\binom{N-s}{k-i}}{\binom{N}{k}}}^M  =\brac{1-\frac{s\binom{N-1}{k-1}}{\binom{N}{k}}}^M\\
&=\brac{1-\frac{ks}{N-s-k+1}\prod_{i=0}^{k-1}\brac{1-\frac{s}{N-i}}}^M.\label{AS}\\
\Pr(\cA_{S,1})
&=\brac{\sum_{i\neq 2}^{k}\frac{\binom{s}{i}\binom{N-s}{k-i}}{\binom{N}{k}}}^M = \brac{1-\frac{\binom{s}{2}\binom{N-2}{k-2}}{\binom{N}{k}}}^M=\\
&= \brac{1-\frac{k(k-1)s(s-1)}{2(N-s-k+2)(N-s-k+1)} \prod_{i=0}^{k-1}\brac{1-\frac{s}{N-i}}}^M.\label{AS1} 
\end{align}
Now we can, for some $r$ (equal to the number of zeroes in $\bsa$), bound the probability of $\cE_{s,\bsa}$ by the product of the RHS of \eqref{AS} with $M$ replaced by $r$ and the RHS of \eqref{AS1} with $M$ replaced by $M-r$. It follows therefore that, after ignoring conditioning on the event $\cB$ that every row of \bA\ contains at least $\g k\s$ ones, we have
$$\Pr(\cE_{s,\bsa})\leq \brac{1-\frac{(1+o(1))k(k-1)s^2e^{-ks/N}}{2(N-s)^2}}^M\leq \brac{1-\frac{k^2 e^{-2k}}{3}}^M\leq (1-e^{-2k})^M.$$
So, in fact, taking account of $\cB$, we have
\beq{BS}{
\Pr(\cE_{S,\bsa}\mid \cB)\leq \frac{\Pr(\cE_{S,\bsa})}{\Pr(\cB)} \leq \frac{(1-e^{-2k})^M}{\Pr(\cB)}.
}
We need a lower bound for $\Pr(\cB)$. By \eqref{rho=l} above, we have,
$$\Pr(\cB)\geb \frac{1}{N^{1/2}}\brac{1-e^{-\l}\sum_{i=0}^{k\g \s-1} \frac{\l^i}{i!}}^N\geq  \frac{1}{N^{1/2}}\brac{1-e^{-(1-\g)^2k\s/3}}^N.$$
Plugging this into \eqref{BS} we see that for large $k$, since $(1-\g)^2\s\geq e^{5k}$,
$$\Pr(\cE_{S,\bsa}\mid \cB)\leq (1-e^{-2k})^M \leq (1-e^{-2k})^{e^{5k}N}.$$
So,
$$\Pr(\exists S,|S|\geq Ne^{-k}:\cE_{S,\bsa})\leq  \sum_{s=Ne^{-k}}^N\binom{N}{s} (1-e^{-2k})^{e^{5k}N}=O(N^{-K}).$$
Finally, if $N/2<|S|\leq N-1$, then the complement $\bar S$ of $S$ is non-empty, and the rows $S$ sum to $\bsa$ if and only if the rows $\bar S$ sum to $\bsb-\bsa$, where $\bsb$ is the row-sum of $\bA$.  But this probability is controlled by the cases above, since $1\leq |\bar S|<N/2$.
\proofend
\subsection{The initial rows of $\bB^{-1}$ have many, but not too many, ones}
We argue next that the rows of $\bB^{-1}$ must contain many ones. Let $\br_1,\br_2,\ldots, \br_{n_2}$ denote the rows of $\bB^{-1}$. We consider its first row $\br_1$. Let $\bb_1,\bb_2,\ldots,\bb_{m_2}$ be the columns of $\bB_1$. Then we must have $\br_1\bb_1=1$ and $\br_1\bb_i=0$ for $i=2,3,\ldots,m_2$. Suppose that $\br_1$ has $s$ ones and let this event be $\cE_0=\cE_0(s)$. Then, for $\cE_0$ to occur there must be $s$ rows of $\bB_1$ whose sum is $(1,0,0,\ldots,0)$.

We apply Lemma \ref{lemrank} to $\bB_1$ with $N=n_2,M=m_2,\g=\frac{n_2}{10m_2}\leq \frac{1}{2}$ and $\bsa=(1,0,0,\ldots,0)$. We consider case (a) and we assume that $s\leq s_0=n_2e^{-k}$. In which case we find, using \eqref{abound}, that
\beq{from}{
\Pr(\cE_0)\leb n^{1/2}\sum_{s=2}^{s_0}\bfrac{se^{2k/3}}{n_2}^{ks/30}=O(n^{-k/50}).
}
Now suppose that $\br_i$ has $\b_in_2$ ones. We can assume from \eqref{from} that 
\beq{bie}{
\b_i\geq \e_0n_2\text{ where }\e_0=e^{-k}.
}
We also need a bound on $1-\b_i$. Again consider $\br_1$. Suppose that this has at least $n_2(1-\e_0)$ ones in positions $S$. Now since each column of $\bB_1$ has exactly $k$ ones, we know that the sum of the rows of $\bB_1$ is either {\bf 0} (if $k$ is even) or {\bf 1}=(1,1,\ldots,1) (if $k$ is odd). (We take care of $k$ even, even though the assumption is still that $k$ is odd.) Thus the $n_2-s$ rows of $\bB_1$ corresponding to $[n_2]\setminus S$ will sum to $(1,0,0,\ldots,0)$ or $(0,1,1,\ldots,1)$ according as $k$ is even or odd. We can apply Lemma \ref{lemrank} once more. This deals with all rows because the probability in \eqref{from} is bounded by $o(n^{-1})$ and so we can use the union bound.
\begin{remark}\label{rem1}
We see that if we fix a positive integer $K$ and if $k$ is sufficiently large, then $\sum_{i\in I}\br_i$ contains at least $s_0$ ones for all $|I|\leq K$. This is because each such $I$ gives us an $\bsa$ with only $|I|$ ones viz. the characteristic vector of $I$. There are $O(n^K)$ such $\bsa$ and the probability bound in \eqref{from} will be small enough to deal with all such $I$ if $K<k/50$.
\end{remark}
\subsection{A few rows of $\bB^{-1}$ are not enough to cover  $[n_2]$}\label{notenough}
Let $S_i,i\in [n_2]$ be the indices of the columns where row $i$ of $\bB^{-1}$ has a one. We will apply the following lemma to the complements of the $S_i$'s. In which case we will have $N=n_2$, $X_i=[n_2]\setminus S_i$ and $\delta=\e_0$.
\begin{lemma}\label{intersection}
Let $X_1,X_2,\ldots,X_N\subseteq [N]$ satisfy $|X_i|\geq \d N$. Let $r$ be a fixed positive integer independent of $N$. If $N$ is sufficiently large, then there exists a set $I\subseteq [N],|I|=r$ and $s=\rdup{\log_2r}$ such that $\card{\bigcap_{i\in I}X_i}\geq \d_s N/2$. Here $\d_0=\d$ and $\d_{i+1}=\d_i^2/4$ for $i\geq 0$.
\end{lemma}
\proofstart
We will assume that $r=2^s$ is a power of two. For general $r$ we take the smallest power of two greater than $r$. This will explain the extra factor of two in the denominator in our lower bound on $\card{\bigcap_{i\in I}X_i}$.

We will prove this by induction on $s$. As a base case, consider $s=1$. Now suppose that for some $t\geq 2$ we find that $|X_t\cap X_i|\leq \d N/(2t)$ for all $i<t$. This implies that $\card{X_t\setminus\bigcup_{i=1}^{t-1}X_i}\geq \d N/2$ and so $\card{\bigcup_{i=1}^tX_i}\geq t\d N/2$. This process must stop after $2/\d$ steps and our induction on $s$ has a base case, i.e. there exists $i,t\leq 2/\d$ such that $|X_i\cap X_t|\geq \d^2N/4$.

Suppose that for some $s$ we can we can find $\set{i_1,i_2,\ldots,i_{2^s}}\subseteq \left[\prod_{i=1}^s (2/\d_i)\right]$ such that $\card{Y_1}\geq \d_s N$ where $Y_1=\bigcap_{j=1}^{2^s}X_{i_j}$. Assuming $N$ is sufficiently large, we can generate a sequence $Y_1,Y_2,\ldots,Y_{2/\d_s}$ where (i) $|Y_i|\geq \d_s N$ for $i=1,2,\ldots,2/\d_s$ and (ii) each $Y_i$ is the intersection of $2^s$ distinct $X_j$ and (iii) no $X_j$ appears in more than one of these intersections. Applying the argument that gave us the base case we see that there exists $i,t\leq 2/\d_s$ such that $|Y_i\cap Y_t|\geq \d_{s+1}N$.
\proofend\\
Putting $X_i=[n_2]\setminus S_i$ for $i\in [n_2]$ we see that we can find for any constant $r$, a set of $r$ rows, such that there are $\Omega(n)$ columns without a one in the union of the rows.  
\subsection{Constructing a representative matrix}\label{largeM}
We now consider the construction of the partition $A_0,A_1,\ldots,A_\ell$ in Section \ref{outline}. Let $R$ denote an arbitrary set of $p$ rows of $\bB^{-1}$. Let $\e_1=2^{-2p}\e_0$, where $\e_0=e^{-k}$, as in \eqref{bie}, and consider the $p\times 2^p$ matrix ${\bD}$ with entries in $\set{0,1}$. The $i$th row $\bu_i$ of $\bD$ is associated with set $S_i$ and the columns of $\bD$ are indexed by ${\boldsymbol \sigma}=(\xi_1,\xi_2,\ldots,\xi_p)\in 2^{[p]}$ and they are associated with an atom $A_{\boldsymbol \sigma}=\bigcap_{j=1}^pS_j^{\xi_j}$ in the Boolean algebra $\cB_R$ generated by the sets $S_i$. Here $\xi_j=\xi_j({\boldsymbol \sigma})=0,1$ and $S_j^1=S_j,S_j^0=\bar{S}_j=[n_2]\setminus S_j$. The columns run over the $2^p$ sequences $\set{0,1}^p$. For each $j\in [n_2]$ there is a unique $\bss=\bss(j)$ such that $j\in A_\bss$ i.e. the $A_\bss$ partition $[n_2]$. Thus $\ell=2^p-1$ here. Further, if $S_{\boldsymbol \sigma}=\bigcap_{i=1}^pS_i^{\xi_i}$ then $S_i$ is partitioned into the parts $S_{\boldsymbol \sigma}$ such that $\xi_i({\boldsymbol \sigma})=1$.

Row $i$ of ${\bD}$ contains a one in position ${\boldsymbol \sigma}$ if $\xi_i({\boldsymbol \sigma})=1$ and $\card{A_{\boldsymbol \sigma}}\geq \e_1n_2$. Otherwise, row $i$ of ${\bD}$ contains a zero in position ${\boldsymbol \sigma}$. We now claim that ${\bD}$ has row rank $p$. 

Fix some $\emptyset\neq I\subseteq [p]$ and let  $\br^I=\sum_{i\in I}\br_i$ and $S_{\oplus}=\set{j:\br^I_j=1}$. Note that Lemma \ref{lemrank} and Remark \ref{rem1} means that we can assume that $|S_\oplus|\geq \e_0n$. Now let $\bos=\sum_{i\in I}\bu_i$ and $S_\bos=\bigcup_{\eta_{\boldsymbol \s}=1}S_{\boldsymbol \sigma}$. We have 
$$|S_\bos|\geq |S_\oplus|-2^p\e_1n_2\geq \e_0n_2-2^p\e_1n_2>0.$$
{\bf Explanation:} When an entry $u_{i,\bss}=1$ this means (among other things) that $j\in S_i$ for all $j\in A_\bss$ and thus $r_{i,j}=1$ for all $j\in A_\bss$. Thus, $S_\bos$ is equal to $S_\oplus$ minus sets of the form $S_{\boldsymbol \sigma}$ where (i) $\xi_i({\boldsymbol \sigma})=1$ for an odd number of $i\in I$ and (ii) $|A_{\boldsymbol \sigma}|\leq \e_1n_2$. 

It follows that there exists ${\boldsymbol \s}$ such that $\bos_{\boldsymbol \sigma}=1$ i.e. $\bos\neq {\bf 0}$. Because $I$ is arbitrary, we see that ${\bD}$ has full row rank.
\subsection{Finishing the proof of Theorem \ref{th1}}\label{finish}
Recall that the minor $M$ can be represented by a $\n\times \m$ matrix $\bR_M$. Let $R$ be a set of row indices where (i) $|R|=\n$ and (ii) $\card{[n]\setminus \bigcup_{i\in R}S_i}\geq \d_s n,s=\rdup{\log_2\n}$ (see Lemma \ref{intersection}). Suppose that $\bc$ is a column of $\bA_m$ not involved in the construction of $\bB$. We say that $\bc$ is a {\em candidate} column if $c_j=0$ whenever $j\in A_{\boldsymbol \sigma}$ for which $|A_{\boldsymbol \sigma}|<\e_1n$. Next let $c_\bss=\sum_{j\in A_\bss}c_j$. If $\bc$ is a candidate column then $\br_i\cdot\bc=\bu_i\cdot\bc_R$ where $\bc_R$ is the column vector with components $c_\bss,\bss\in 2^{[\n]}$. (Remember that $\br_i$ is row $i$ of $\bB^{-1}$ and that $\bu_i$ is row $i$ of ${\bD}$.) 
For a column \bx\ of $\bA_m$, let $\f_R(\bx)$ be the column \bx\ restricted to the $\n$ rows of $R$. Let $\bc_1$ be the first column of the target matrix $M$ and let $\bc$ be a random candidate column.  Let $\bv$ satisfy ${\bD}\bv=\bm_1$ and $v_\bss=0$ if $|A_{\boldsymbol \sigma}|<\e_1n$. Assume also that $\bv$ has at most $\n$ ones and that $k\geq \n$. There are always such solutions. Then we have
\beq{prob}{
\Pr(\f_R(\bB^{-1}\bc)=\bm_1)=\Pr(\br_i\cdot\bc=\bu_i\cdot\bc_R=m_{1,i},i\in R)\geq\Pr(\bc_R=\bv)\geq\e_1^k.
}
{\bf Explanation of second inequality:} Let $J=\set{\bss: v_{\boldsymbol \sigma}=1}$. Each index ${\boldsymbol \sigma}$ corresponds to a set $A_{\boldsymbol \sigma}$ of size at least $\e_1n$. Now we will have $\bc_R=\bv$ if column $\bc$ has a single one in each $A_{\boldsymbol \sigma},{\boldsymbol \sigma}\in J$ and its remaining ones $A_0=[n_2]\setminus \bigcup_{i\in R}S_i$. All of the sets where we need to place ones are of size at least $\e_1n_2$ and \eqref{prob} follows.

It follows from this that we can find a copy of $M$ w.h.p. by examining a further $\om$ random columns, where $\om=\om(n)\to\infty$, is arbitrary. This completes the proof of Theorem \ref{th1}.
\section{$k$ even}\label{keven}
We now examine the adjustments needed for the case of $k$ even. The problem here is that the rows of $\bA_m$ now sum to zero and so we cannot construct $\bB$ in quite the same way as for $k$ odd. 

Going back to Section \ref{extendB} we define $\bL_1,\bL_2,\bL_3$ in the same way, but now we can only say that w.h.p. the rank of $\bL_1$ is $n_3^*=n_3-1$. So now we choose $i\in I_2$ such that if the matrix $\bL_2^*=[\bB_1^*:\bL_3^*]$ is obtained from $\bL_2$ by deleting row $i$ then $\bL_3^*$ has full row rank. We claim now that w.h.p. $\bL_2^*$ also has full row rank. Suppose now that there exists $J\subseteq I_1\setminus\set{i}$ such that $\sum_{j\in J}\ba_j=0$. Then we must have $J\cap (I_2\setminus \set{i})=\emptyset$. For otherwise, $\bL_3^*$ does not have full row rank. But then $J\subseteq I_1\setminus I_2$ and by \eqref{C1==0} we can assume that $|J|\leq n_2e^{-2k}$ and then we can apply Lemma \ref{lemrank}(a) to $J$ and $\bB_1$ to get a contradiction w.h.p.

Then we let $\bB^*$ be obtained from $\bB_1$ by removing row $i$ and then we can extend it to an $n_2^*\times n_2^*$ non-singular submatrix of $\bL_2^*$. We need to argue that $\bB^*$ has many ones and zeros in each row. For this we need to argue about sums of rows of the matrix $\bB_1^*$ which has $k$ or $k-1$ ones in each column and at least $k/10$ ones in each row. The row we deleted from $\bB_1$ came by considering $\bL_1$ which is independent of $\bB_1$ and so the ones in each column of $\bB_1^*$ are still randomly chosen subject to the row constraints. We can then argue via Lemma \ref{lemrank} that ${(\bB^*)}^{-1}$ has many ones and zeros in each row.

We write
$$\bA_m=\left[\begin{array}{cccc}
\bB^*&\bL_1^*&\bC^*&\bR^*\\
\bu_1&\bu_2&0&\bu_3
\end{array}\right]$$
where $[\bu_1,\bu_2,0,\bu_3]$ is row $i$ and $\bC^*$ comprises the unviewed random columns that appear where there is a zero in row $i$. $\bR^*,\bu_3$ comprise the rest of the matrix. Now let $\widehat \bB$ be the $n_2\times n_2$ matrix obtained from $\bB^*$ by adding a column $\be_{n_2}$ and a row $\be_{n_3}^T$ where $\be_{n_2}$ has a unique one in position $n_3$. 

The number of ones in a row of $\bX$ is dominated by the binomial $Bin(n/4,k/n)$ and so w.h.p. the maximum number of ones in any row is $O(\log n)$. Then we write
$${\widehat \bB}^{-1}\bA_m=\left[\begin{array}{cccccc}
\bI_1&0&\widehat \bL_{1,1}&\widehat \bC_{1,1}&\widehat \bC_{1,2}&\widehat \bR_1\\
0&\bI_2&\widehat \bL_{2,1}&0&\widehat \bC_{1,2}&\widehat \bR_2\\
0&\bu_{1,2}&\bu_2&0&0&\bu_3
\end{array}\right].$$
Here we have split $\bu_1$ into $[\bu_{1,1},\bu_{1,2}]$ where $\bu_{1,1}=0$ and $\bu_{1,2}$ is an all ones vector of dimension $O(\log n)$. And then the matroid $\cM$ associated with $\bA_m$ has a minor isomorphic to $M$ if $\bR_M$ appears in $\widehat \bC_{1,1}$. The argument for this is covered by the case $k$ odd, concentrating on the sub-matrix $[\bI_1:\widehat \bC_{1,1}]$.
\section{Further Questions}
We have shown that $\bA_m$ contains a copy of an arbitrary fixed binary matroid as a minor under the assumption that $k,m/n$ are sufficiently large. It would be of interest to reduce $k$, perhaps to three, and to get precise estimates for the number of columns needed for some fixed matroid, the Fano plane for example. In this way we could perhaps get the precise number of columns needed to make the random matroid associated with $\bA_m$, non-graphic or non-regular, w.h.p. Behavior of random matroids over fields other than $\GF_2$ are also an interesting target.

{\bf Acknowledgement:} We thank Peter Nelson for for helpful discussions.

\appendix
\section{Proof of \eqref{rho=l}}
Let $\brho$ be the vector of row counts in $\bX$ and let $A,B$ be arbitrary positive integers,
$$S = \Bigl\{ \brho \in [M]^N \,\Big|\sum_{1\leq j \leq N} \r_j = A \mbox{ and }
\forall j,\, \r_j \geq B\Bigr\}.$$   
Fix $\vec \xi \in S$.   Then, if $\Pr_1$ refers to a random choice from $S$,
$$\Pr(\brho = \vec \xi) =
\left( \frac{M!}{\xi_1! \xi_2! \ldots \xi_N! }\right)
\bigg/
\left( \sum_{\brho\in S} \frac{M!}{\r_1! \r_2! \ldots \r_N!} \right).$$
On the other hand, if $\Pr_2$ refers to a random choice via independent Poisson,
\begin{align*}
\Pr_2\left(\brho = \vec \xi \; \bigg|\; \sum_{1\leq j \leq N} \r_j =  A\right)
&=\left( \frac{\prod_{1\leq j\leq N}\l^{\xi_j}}{f_B(\l)\xi_j!}\right)\bigg/\left( \sum_{\brho\in S} \prod_{1\leq j\leq N}\frac{\l^{\r_j}}{f_B(\l) \r_j!}\right)\\
&=\left(   \frac{ f_B(\l)^{-N} \l^{s}}{\xi_1! \xi_2! \ldots \xi_N!} \right)\bigg/\left( \sum_{\brho\in S} \frac{f_B(\l)^{-N} \l^{s}}{\r_1! \r_2! \ldots \r_N!} \right)\\
&=\Pr_1(\brho=\vec \xi).
\end{align*}
\proofend

\begin{thebibliography}{99}
\bibitem{AY} J. Altschuler and E. Yang, Inclusion of Forbidden Minors in Random Representable Matroids, arXiv:1507.05332 [math.CO].

\bibitem{AFP} J. Aronson, A.M. Frieze and B. Pittel, Maximum matchings in sparse random graphs: Karp-Sipser revisited, {\em Random Structures and Algorithms} 12 (1998) 111-178.

\bibitem{B} B. Bollob\'as, Random Graphs, First Edition, Academic Press, London 1985, Second Edition, Cambridge University Press, 2001.

\bibitem{BPP} N. Bansal, R.A. Pendavingh and J.G. van der Pol, On the number of matroids, {\em Combinatorica} 35 (215) 253-277.

\bibitem{DM} O. Dubois and J. Mandler, The 3-XORSAT Threshold, {\em Proceedings of the 43rd Annual IEEE Symposium on Foundations of Computer Science} (2002) 769-778.

\bibitem{Co1} C. Cooper, On the rank of random matrices, {\em Random Structures and Algorithms} 16 (2000) 209-232. 

\bibitem{CDCCores} C. Cooper, The cores of random hypergraphs with a given degree sequence, {\em Random Structures and Algorithms} 25 (2004) 353-375.

\bibitem{FK} A.M. Frieze and M. Karo\'nski, Introduction to Random Graphs, Cambridge University Press 2015.

\bibitem{growth} J. Geelen, J.P.S.~Kung, and G.~Whittle, Growth rates of minor-closed classes of matroids, {\em Journal of Combinatorial Theory, Series B} {\bf 99} (2009) 420-427.
 
\bibitem{JLR} S. Janson, T. {\L}uczak and A. Ruci\'nski, Random Graphs, John Wiley and Sons, New York, 2000.

\bibitem{K15} M. Kahle, Topology of random simplicial complexes: a survey, {\em AMS Contemporary Volumes in Mathematics}, 

\bibitem{KO} D. Kelly and J. Oxley, On random representable matroids, {\em Studies in Applied Mathematics} 71 (1984) 181-205.

\bibitem{Kung} J.P.S.~Kung, The long-line graph of a combinatorial geometry. II. Geometries representable over two fields of different characteristics, {\em Journal of Combinatorial Theory, Series B} {\bf 50} (1990) 41--53.
  
\bibitem{Mo} M. Molloy,  Cores in random hypergraphs and Boolean formulas, {\em Random Structures and Algorithms} 27 124-135 (2005). 

\bibitem{LM06} N. Linial and R. Meshulam, Homological connectivity of random 2-complexes, {\em Combinatorica} 26 (2006) 475-487.

\bibitem{Ox} J. Oxley, Matroid Theory, Second Edition, Oxford University Press, New York, 2011.

\bibitem{OSWW} J. Oxley, L. Lowrance, C. Semple and D. Welsh, On properties of almost all matroids, {\em Advances in Applied Mathematics} 50 (2013) 115-124.

\bibitem{M} M. Molloy, Cores in random hypergraphs and random formulas, {\em Random Structures and Algorithms} 27 (2005) 124-135.

\bibitem{Nel} P. Nelson,  Almost all matroids are nonrepresentable, {\em Bulletin of the London Mathematical Society} 50 (2018) 245-248.

\bibitem{PP} R. Pendavingh and J.  van der Pol,  On the number of bases of almost
all matroids, arXiv preprint arXiv:1602.04763.

\bibitem{Tu} W.T. Tutte, A homotopy theorem for matroids. I, II, {\em Transactions of the American Mathematical Society} 88 (1958) 144–174.

\bibitem{Welsh} D. Welsh, Matroid Theory, Academic Press, 1976.

\end{thebibliography}
\end{document}